\documentclass[12pt]{article}

\usepackage{amsmath, amsfonts, amssymb}
\usepackage{latexsym}
\usepackage{graphicx}
\usepackage{color}
\usepackage{url}

\newtheorem{thm}{Theorem}[section]

 \newtheorem{prop}[thm]{Proposition}

\newcommand{\bx}{{\mathbf x}}
\newcommand{\by}{{\mathbf y}}

\newcommand{\bob}{{\mathbf b}}

\newcommand{\ba}{{\mathbf a}}
\newcommand{\bo}{{\mathbf 0}}

\newcommand{\bZ}{{\mathbf Z}}

\newcommand{\cI}{{\mathcal I}}

\newcommand{\al}{\alpha}
\newcommand{\be}{\beta}

\renewcommand{\thefootnote}{\fnsymbol{footnote}}

\begin{document}
\begin{center}
{\Large \bf Representing integers by multilinear polynomials}

\bigskip
{\bf Albrecht B\"ottcher and Lenny Fukshansky}
\end{center}

\begin{quote}
\footnotesize{Let $F(\boldsymbol x)$ be a homogeneous polynomial in $n \ge 1$ variables of
degree $1 \leq d \leq n$ with integer coefficients so that its degree in every variable is equal to $1$.
We give some sufficient conditions on $F$ to ensure that for every integer $b$ there exists
an integer vector $\boldsymbol a$  such that $F(\boldsymbol a) = b$. The conditions
provided also guarantee that the vector $\boldsymbol a$ can be found in a finite number of
steps.
}
\end{quote}

\let\thefootnote\relax\footnote{\hspace*{-7.5mm} MSC 2010: Primary 11D85, 11C08, 11C20; Secondary 11G50}
\let\thefootnote\relax\footnote{\hspace*{-7.5mm} Keywords: polynomials, integer representations, unimodular matrices, linear and multilinear forms}
\let\thefootnote\relax\footnote{\hspace*{-7.5mm}  Fukshansky acknowledges support by Simons Foundation grant \#519058.}

\section{Introduction and main result}

The problem of finding solutions to a given polynomial equation with integer coefficients goes back to the work
of Diophantus, resulting in these equations being called the Diophantine equations. While the linear case likely dates back to Diophantus himself, the quadratic equations have been systematically studied by Gauss and famously led to his composition law for the binary quadratic forms. More generally, Hilbert's 10th problem, in its contemporary formulation, asks whether there exists an algorithm to decide if
a given polynomial equation with integer coefficients has a (nontrivial) integer solution.
Matiyasevich's famous theorem~\cite{mat} of 1970 (building on previous work of others)
gives a negative answer to this problem. In fact, J. P. Jones~\cite{jp_jones} proved in 1980 that the question
whether a general Diophantine equation of degree four or larger has a solution in positive integers is already undecidable,
and not much else is known for polynomials of degree~$\ge 4$. 

One possible approach to Hilbert's 10th problem is through {\it search bounds}. Suppose we can prove that if a polynomial
$F(x_1,\dots,x_n)$ with integer coefficients has an integer zero~$\ba \in \bZ^n \setminus \{\bo\}$, then it has one with
$|\ba| \le C(n,F)$,
where $|\cdot |$ stands for, say, the sup-norm and $C(n,F)$ is some explicitly given function depending on $n$ and $F$.
Since the set of integer points $\ba \in \bZ^n$ satisfying this condition is finite, we can simply search through
all of them checking whether any one of these vectors is a zero of $F(x_1,\dots,x_n)$. This approach will either
produce a solution or prove that one does not exist. A survey of known results on search bounds can be found in~\cite{masser}. While Matiyasevich's theorem guarantees that search bounds
are not possible in general, an overview of results of this type in the quadratic case can be
found in~\cite{cassels_overview}, and the current state of the art for the cubic case is in~\cite{bde}, as well as the references therein.

In this note we study a special class of polynomials of arbitrary degree. Let $n \geq 1$ be an integer and let us define
$[n] := \{1,\ldots,n\}$.
Given an integer $d$ with $1 \leq d \leq n$, we put
$\cI_d(n) := \left\{ I \subseteq [n] : |I| = d \right\}$.
For each indexing set $I = \{i_1,\dots,i_d\} \in \cI_d(n)$ with $1 \le i_1 < \dots < i_d \le n$, we define the monomial
$x_I$ in the variables $x_{i_1},\dots,x_{i_d}$ out of $x_1,\dots,x_n$ as $x_I := x_{i_1} \cdots x_{i_d}$.
An {\em integer multilinear $(n,d)$-form} is a polynomial of the form
$$F(x_1,\dots,x_n) = \sum_{I \in \cI_d(n)} f_I x_I,$$
where the coefficients $f_I$ are integers for all $I \in \cI_d(n)$. Such an $F$ is a homogeneous polynomial
in $n$ variables of degree $d$ which has degree $1$ in each of the variables $x_1,\dots,x_n$. We will say that
$F$ {\em represents} an integer $b$ it there exists an integer vector $\ba \in \bZ^n$ such that
$F(\ba) = b$.
Under what conditions on $F$ does such a polynomial represent all integers? The first observation is that the coefficients
$f_I$ of $F$ must be relatively prime: if
$g = \gcd(f_I)_{I \in \cI_d(n)} > 1$,
then $g$ must divide $F(\ba)$ for every $\ba \in \bZ^n$, and hence an integer $b$ that is not a multiple of $g$ is not represented by~$F$.
We will say that our form {\em is coprime} if $\gcd(f_I)_{I \in \cI_d(n)} = 1$.

We will provide some sufficient conditions for a multilinear $(n,d)$-form $F$ to represent all integers.
Further, our results are effective in the sense that we provide algorithms that yield an integer solution $\ba$ of the equation
$F(\ba)=b$ (theoretically but not necessarily practically) in a finite number of steps.
Here is our main result.

\begin{thm} \label{main1} Let $F(\bx)$ be a coprime integer multilinear $(n,d)$-form.
Suppose in addition that at least one of the following two conditions holds:

$\quad$ {\rm (a)} The nonzero coefficients of $F$ are pairwise coprime,

$\quad$ {\rm (b)} $n = d+1$ and $F$ has a pair of coprime coefficients.

\noindent
Then $F$ represents all integers. Further, for each $b \in \bZ$ there exists an $\ba \in \bZ^n$ such that $F(\ba) = b$ and
$$|\ba|  \le |b| \left( 2|F| \right)^{ d!\,{\rm e}},$$
where $|\ba|=\max_{1 \le i \le n}|a_i|$, $|F| = \max_{I \in \cI_d(n)} |f_I|$, and ${\rm e} = 2.71828\ldots$.
\end{thm}

We will prove the theorem with an exponent $\nu_d$ that is slightly  sharper than $d! {\rm e}$. We neither know what the sharpest exponent could be like nor do we have a lower bound for it.
We also remark that neither (a) nor (b) is a necessary condition. For example, it is well known (see~\cite{oleary} and the references therein) that a coprime integer linear $(n,1)$-form
$F(\bx)=f_1x_1+ \cdots+f_nx_n$ represents all integers even if no pair of coefficients is coprime, and there are $\ba \in \bZ^n$ such that $F(\ba)=b$
and $|\ba| \le |b|$. To give another example,
the coprime $(3,2)$-form $F(x,y,z)=6xy+10xz+15yz$ represents every integer
although the $\gcd$ of each pair of coefficients is greater than $1$. This is the case $p=5$ of the following observation.

\begin{prop} \label{main2}
If $p \ge 5$ is an integer and neither $2$ nor $3$ divides $p$, then the polynomial $F(x,y,z)=6xy+2pxz+3pyz$ represents every integer.
\end{prop}

Next we turn to a special class of multilinear forms, one for which it is easy to establish an ``if and only if" result. 
The following observation is closely related to Lemma~2 on page~15 of~\cite{cassels}, 
which is also recorded in~\cite{zurich} (Lemma~2), where it is referred to as the PID version of Quillen-Suslin's theorem. 
We derive it from a more general result by Xingzhi Zhan~\cite{zhan} (Theorem~1).

\begin{thm} \label{main3}
Let $A=(a_{ij})\in \bZ^{n \times s}$ with $1 \le s < n$ and consider the integer multilinear
$(n(n-s),n-s)$-form
\[F(\by)={\rm det}\left(\begin{array}{llllll}
a_{11} & \ldots & a_{1s} & y_{11} & \ldots & y_{1(n-s)}\\
\vdots & & \vdots & \vdots & & \vdots\\
a_{n1} & \ldots & a_{ns} & y_{n1} & \ldots & y_{n(n-s)}
\end{array}\right).\]
This polynomial represents all integers if and only if
the minors of order $s$ of $A$ are coprime. If this is the case,
there is a finite algorithm to find an integer solution of $F(\by)=b$
for all $b \in \bZ$ simultaneously.
\end{thm}

The algorithm mentioned in the previous theorem is basically the algorithm for constructing
the Smith normal form $A=USV$ of the matrix $A$; here $U$ and $V$ are unimodular integer
matrices, that is, integer matrices with determinant $\pm 1$, while $S$ is the diagonal matrix
given by the invariant factors of $A$. We can also produce explicit search bounds as follows.
Suppose the minors of order $s$ of $A$ are coprime and denote by $D$ the minimum of the absolute values
of the nonzero minors. Let $\al$ be the maximum of the absolute values of the entries of $A$
and put $\be =(n-1)!D^n+1$. Then there is a $\by \in \bZ^{n(n-s)}$ such that $F(\by)=b$
and
$$|\by| \le n^2 |b| \al \be (\be+1)^{n-2}.$$
This bound can be obtained by employing careful basis change arguments as on pages 10 -- 15
of Cassels' book~\cite{cassels}. However, as the algorithms for getting the Smith normal
form of $A$ are more efficient than the search based on this bound, we will omit the proof
of this bound here.

We end this introductory section with an open problem.

\medskip
{\bf Question.} {\em Does there exist a coprime integer multilinear $(n,d)$-form that does not represent all integers?}

\section{Proofs and further results}

{\bf Proof of Theorem \ref{main1}(a).} By the remark after Theorem~\ref{main1}, we may assume that $d \ge 2$. We define
\begin{equation}
\label{mu}
\nu_{d} = \sum_{k=0}^{d} \frac{d!}{k!}
\end{equation}
and will show the theorem with the bound
\begin{equation}
\label{abnd}
|\ba| \le |b| \left( 2|F| \right)^{\nu_{d}}.
\end{equation}
As $\nu_d <d!\,{\rm e}$, this is actually sharper than the bound given in Theorem~\ref{main1}.

Since $F(\bx) := F(x_1,\dots,x_n)$ is homogeneous, $F(\bo) = 0$. Hence from here on we assume that
$b \neq 0$. First suppose that $F(\bx)$ has only one monomial, i.e.,
$$F(\bx) = f_I \prod_{i \in I} x_i$$
for some $I \subseteq [n]$ and $f_I \in \bZ$. Since the $\gcd$ of the coefficients of $F$ is $1$,
we must have $f_I = \pm 1$. Take some $j \in I$ and put $a_j = \pm b$ and $a_i = 1$ for $i \in I \setminus\{j\}$.
We so obtain a vector $\ba \in \bZ^n$ such that $F(\ba) = b$ and
$|\ba| = \max\{1, |b|\} = |b|$, which is smaller than the bound \eqref{abnd}.

Next assume that $F(\bx)$ has exactly two monomials, i.e.,
$$F(\bx) = f_{I_1} \prod_{i \in I_1} x_i + f_{I_2} \prod_{i \in I_2} x_i$$
for some $I_1,I_2 \subseteq [n]$ and coprime $f_{I_1}, f_{I_2} \in \bZ$. Then the index sets $I_1$ and $I_2$ must be
distinct (since otherwise there would be only one monomial) of the same cardinality $d$, and
so there must exist some $k \in I_1 \setminus I_2$ and $m \in I_2 \setminus I_1$. Let $a'_k, a'_m \in \bZ$ be such that
$$a'_k f_{I_1} + a'_m f_{I_2} = 1.$$
The Euclidean algorithm allows us to find such $a'_k, a'_m$ with
$$|a'_k|, |a'_m| \le \max \{ |f_{I_1}|, |f_{I_2}| \}.$$
Letting $a_k = b a'_k$, $a_m = b a'_m$, and $a_i=1$ for $i \neq k,m$, we get
$$F(\ba) = a_k f_{I_1} + a_m f_{I_2} = b$$
with $|\ba| \leq |b| |F|$, which is again smaller than the bound \eqref{abnd}.

We now argue by induction on $\ell \ge 1$, the number of monomials of $F$. Since the base of induction is already established,
we assume that $\ell \geq 3$ and that the result is proved for polynomials with no more than $\ell - 1$ monomials.
First notice that we can assume without loss of generality that $F$ depends on all variables (if not, then $F$ is a polynomial in $< n$ variables)
and that no variable is present in all monomials (if it is, then just set it equal to $1$).
Let $d \ge 2$ be the degree of $F$. Every monomial is indexed by a subset $I$ of $[n]=\{1,\ldots,n\}$ of cardinality $d$.

Suppose first that the variable $x_1$ is present in $\ell-1$ monomials. We then may write
\begin{equation}
F(\bx) = x_1 G(x_2,\dots,x_n) + f_I \prod_{i \in I} x_i, \label{FG}
\end{equation}
where $I \subset \{ 2,\dots,n \}$ with $|I|=d$ and $G$ is a homogeneous polynomial of degree $d-1$
that is linear in each of the $n-1$ variables with pairwise coprime integer coefficients.
By the induction hypothesis, there exists a vector $\ba' =(a_2, \ldots, a_n) \in \bZ^{n-1}$ such that $G(\ba') = 1$ and
$$|\ba'| \le |1|(2 |G|)^{\nu_{d-1}} \le (2 |F|)^{\nu_{d-1}}.$$
Put $a_1=b-f_I \prod_{i \in I}a_i$. Then
\begin{equation}
F (a_1, \ba')
= \left( b - f_I \prod_{i \in I} a_i \right) G(\ba') + f_I \prod_{i \in I} a_i = b,\label{FGb}
\end{equation}
that is,  $F(\ba) = b$ for $\ba=(a_1,a_2, \ldots, a_n)$. Furthermore,
$$|\ba| \le |b| + |f_I| |\ba'|^{d} \le 2 |b| |F| |\ba'|^{d}$$
since $|b|$, $|f_I| |\ba'|^{d}$ are positive integers and $f_I$ is a coefficient of~$F$. Therefore
$$|\ba| \le 2 |b| |F| \left( 2 |F| \right)^{\nu_{d-1}d} = (2 |F|)^{1 + d \nu_{d-1}}\ |b|,$$
and because, by (\ref{mu}),
$$1 + d \nu_{d-1} = 1 + d\sum_{k=0}^{d-1} \frac{(d-1)!}{k!} = \sum_{k=0}^{d} \frac{d!}{k!}=\nu_{d},$$
we obtain the bound~\eqref{abnd}.

On the other hand, assume that $x_1$ is not present in at least two different monomials.
Then set $x_1=0$ and apply the induction hypothesis to the resulting polynomial
\begin{equation}
P(x_2,\ldots,x_n) := F(0, x_2, \ldots,x_n) \label{PG}
\end{equation}
in $n-1$ variables. This polynomial has no more than $\ell -1$ and no fewer than two monomials
and satisfies all the other conditions of the theorem. Take $\ba' \in \bZ^{n-1}$ to be the point
guaranteed by the induction hypothesis, so that $P(\ba') = b$ and
\begin{equation}
\label{P_bnd}
|\ba'| \leq (2 |P|)^{\nu_{d}}\ |b| \le (2 |F|)^{\nu_{d}}\ |b|.
\end{equation}
Setting $\ba$ to be $\ba'$ with inserted $0$ in the first coordinate, we obtain the necessary
solution~$F(\ba) = b$ with $|\ba| = |\ba'|$ bounded as in~\eqref{P_bnd}, which gives the bound~\eqref{abnd}.
$\;\:\square$

\medskip
{\bf Proof of Theorem \ref{main1}(b).}
We argue by induction on $d \geq 1$. As said, if $d=1$, then $n=2$ and
$F(x_1,x_2) = f_1 x_1 + f_2 x_2$
with $\gcd(f_1,f_2) = 1$. Thus, the result follows from the Euclidean algorithm.

Suppose now $d \ge 2$. Since $n = d+1 \ge 3$, the set $\cI_d(n)$ consists of the indexing sets
$I(k)=[n]\setminus\{k\}$ with $1 \le k \le n$, and so
$$F(x_1,\dots,x_n) = \sum_{k=1}^n f_{I(k)} x_{I(k)}.$$
Since $F$ has a pair of coprime coefficients, there must exist $1 \leq j < m \leq n$
such that $\gcd(f_{I(j)}, f_{I(m)}) = 1$. Assume without loss of generality that $j=n-1$, $m=n$,
and notice that each monomial $x_{I(k)}$ for $k \neq 1$ is divisible by $x_1$. Thus, writing $I'(k) = I(k) \setminus \{ 1 \}$
we obtain
$$F(x_1,\dots,x_n) = x_1 G(x_2, \ldots, x_n) + f_{I(1)} x_{I(1)}
= x_1 G(x_2, \ldots, x_n)+ f_{I(1)} \prod_{i=2}^n x_i$$
with
$$ G(x_2, \ldots, x_n)=\sum_{k=2}^{n} f_{I(k)} x_{I'(k)}.$$
The polynomial $G$
is a coprime integer multilinear $(n-1,d-1)$-form with $n-1=(d-1)+1$ and
$G$ still has the same pair of coprime coefficients $f_{I(n-1)}, f_{I(n)}$.
We can therefore apply the induction hypothesis to $G$ and can argue in the same way as in the proof
of Theorem~\ref{main1}(a) to get the desired result. $\;\:\square$

\medskip
{\bf Proof of Proposition \ref{main2}.}
For $z \in \{\pm 1\}$, the equation $6xy+2pxz+3pyz=b$ is equivalent to the equation
\[(2x+pz)(3y+pz)=b+p^2.\]
Let $b+p^2=2^{\al}m$ with an integer $\al \ge 0$ and an odd integer $m$.
If $p \equiv 1 \mod 3$, then $3y+p=2^\al$ has an integer solution $y_1$ for $\al$ even and $3y-p=2^\al$ has an integer solution
$y_2$ for $\al$ odd. The equations $2x\pm p=m$ always have an integer solution $x_0$. It follows that
$F(x_0,y_1,1)=b$ for $\al$ even and $F(x_0,y_2,-1)=b$ for $\al$ odd.
If $p \equiv -1 \mod 3$, then $3y+p=2^\al$ has an integer solution $y_1$ for $\al$ odd and $3y-p=2^\al$ has an integer solution
$y_2$ for $\al$ even. The equations $2x\pm p=m$ again have an integer solution $x_0$. It follows that
$F(x_0,y_1,1)=b$ for $\al$ odd and $F(x_0,y_2,-1)=b$ for $\al$ even. $\;\:\square$

\medskip
In the specific case when $F$ is a quadratic form, the investigation of bounds for the size of solutions to
equations like $F(\bx) = b$ has a long history: see~\cite{cassels_overview} for a survey of this area.
In particular, Theorem~1 of~\cite{dietmann} gives a bound of the size $|\ba|$ of a solution vector~$\ba$
in case the quadratic form $F$ is nonsingular: the exponent on $|F|$ and $|b|$ in that bound is linear in $n$.
For $d=2$, the bound~(\ref{abnd}) becomes $|\ba| \le |b|(2|F|)^{\nu_2}=|b|(2|F|)^5$.
This bound can be slightly improved.

\begin{prop} \label{main4} Let $F$ be a quadratic form in $n \geq 2$ variables which is linear
in each of the variables, i.e., it contains no diagonal terms. Assume also that the coefficients
of $F$ are pairwise coprime. Then for every $b \in \bZ$ there exists an $\ba \in \bZ^n$ such that $F(\ba) = b$ and
$$|\ba| \le |b|+ |F|^3.$$
\end{prop}

{\em Proof.} Suppose $F$ contains exactly $\ell \ge 2$ monomials and $x_1$ occurs in exactly $\ell-1$ monomials.
We then have~(\ref{FG}) with a linear form $G$ and a set $I \subseteq \{2, \ldots,n\}$ of cardinality $2$.
There is a vector $\ba'=(a_2, \ldots, a_n) \in \bZ^{n-1}$ such that $G(\ba')=1$ and $|\ba'| \le |G| \le |F|$.
With $a_1=b-f_I \prod_{i \in I}a_i$ we get $F(a_1,\ba')=b$ as in~(\ref{FGb}) and clearly,
$$|(a_1, \ba')| \le |b|+|f_I| |\ba'|^2 \le |b|+|F|\,|F|^2=|b|+|F|^3.$$
If there are two monomials in which $x_1$ is not present, we put $x_1=0$ and repeat the argument for the
polynomial~(\ref{PG}). $\;\:\square$

\medskip
Theorem \ref{main3} is the $r=n$ case of the following more general result.

\begin{thm} \label{main5} Let $A=(a_{ij})\in \bZ^{r \times s}$ with $1 \le s \le r \le n$ and $s <n$. Consider
the polynomial
\[F(\bx,\by)={\rm det}\left(\begin{array}{llllll}
a_{11} & \ldots & a_{1s} & y_{11} & \ldots & y_{1(n-s)}\\
\vdots & & \vdots & \vdots & & \vdots\\
a_{r1} & \ldots & a_{rs} & y_{r1} & \ldots & y_{r(n-s)}\\
x_{11} & \ldots & x_{1s} & y_{(r+1)1} & \ldots & y_{(r+1)(n-s)}\\
\vdots & & \vdots & \vdots & & \vdots\\
x_{(n-r)1} & \ldots & x_{(n-r)s} & y_{n1} & \ldots & y_{n(n-s)}
\end{array}\right).\]
This polynomial represents all integers if and only if one of the following two conditions
is satisfied:

\smallskip
{\rm (i)} $r+s \le n$,

\smallskip
{\rm (ii)} $r+s >n$ and the minors of order $r+s-n$ of $A$ are coprime.

\smallskip
\noindent
If {\rm (i)} or {\rm (ii)} holds, there is a finite algorithm to find an integer solution of $F(\bx,\by)=m$
for all $b \in \bZ$ simultaneously.
\end{thm}

{\em Proof.} Zhan \cite{zhan} showed that the equation $F(\bx,\by)=1$ has an integer solution if and only if either (i) is satisfied or if (ii)
holds and $A$ has at least $r+s-n$ invariant factors equal to $1$. Let $k:=r+s-n \ge 1$. For $j=1, \ldots, s$, denote by $d_j$ the
greatest common divisor of the $j \times j$
minors of $a$. Put $d_0=1$ and let $q={\rm rank}\,A$. We then have
\[d_0\,|\, d_1 \, | \ldots \, | d_q,\]
and the invariant factors of $A$ are
\[s_1=\frac{d_1}{d_0}, \quad s_2=\frac{d_2}{d_1}, \quad \ldots, s_q=\frac{d_q}{d_{q-1}}.\]
If the minors of order $k$ are coprime, then
then $d_k=1$, so $d_1= \ldots = d_{k-1}=1$,  and hence $q \ge k$ and $s_1= \ldots =s_k =1$,
that is, $A$ has $k$ invariant factors equal to $1$.
Conversely, suppose $k$ invariant factors of $A$ equal $1$. Since also
\[s_1\,|\, s_2 \, | \ldots \, | s_q,\]
it follows that $s_1 = \ldots =s_k =1$, implying that
$q \ge k$ and $d_k=1$, i.e., the minors of order $k$ are coprime.
Thus, if $k=r+s-n \ge 1$, then the coprimeness of the $k \times k$
minors is equivalent to the existence of $k$ invariant factors of magnitude $1$.

Zhan \cite{zhan} constructed integer vectors $\bx, \by$ with $F(\bx,\by)=1$ explicitly in terms of the Smith normal form $A=USV$.
As the Smith normal form can be computed with finitely many steps in integer arithmetic, an integer solution of  $F(\bx,\by)=1$ can
be found in this way. Having a solution with  $F(\bx,\by)=1$, replacement of the first column $\by_1$ of $\by$ by $b\by_1$
yields an integer solution of  $F(\bx,\by)=b$. $\;\:\square$

\medskip
Here is one more simple observation about polynomials representing all integers and being related to unimodular matrices.

\begin{thm}
Let $L_1(x), \ldots, L_m(x) \in \bZ[x_1, \ldots, x_n]$ be integer linear forms in $n$ variables and let $n > m$.
Denote by $A \in \bZ^{m \times n}$ the matrix whose rows are the coefficient vectors of $L_1, \ldots, L_m$.
Suppose the $m \times m$ minors of $A$ are coprime. Then the polynomial
$$ F(\bx)=\prod_{i=1}^m L_i(\bx)$$
represents all integers. Furthermore, for $b \in \bZ$, put $\bob=(b,1,\ldots,1)^\top \in \bZ^m$
and define $\mu(A,b)$ as the maximum of the absolute values of the $m \times m$ minors of the augmented matrix $(A\: \bob)$.
Then there exists an $\ba \in \bZ^n$ such that $F(\ba)=b$ and $|\ba| \le \mu(A,b)$.
\end{thm}

{\em Proof.}  Define
$\gcd(A)$ and $\gcd(A,b)$ to be the greatest common divisors of the $m \times m$ minors of $A$
and $(A \:\bob)$, respectively. Consider the linear system $A\bx = \bob$. A theorem of Ignaz Heger~\cite{heger}
(see also~\cite{oleary}) states that this system has an integer solution if and only
if $\gcd(A) = \gcd(A, b)$. This is clearly the case if the $m \times m$ minors of $A$ are coprime: in this case
$\gcd(A) = 1$, and so $1 \le \gcd(A,b) \le \gcd(A)=1$. The main theorem of~\cite{borosh} therefore guarantees
the existence of a solution vector $\ba \in \bZ^m$
to this system with $|\ba| \le \mu(A,b)$. It follows that
$$F(\ba)=\prod_{i=1}^m L_i(\ba)=b \times 1 \times \cdots \times 1 = b,$$
as desired. $\;\:\square$

\medskip
We remark that if the $m \times m$ minors of $A$ are not coprime, then the polynomial $F(\bx)$ may not be
representing all integers. Indeed, consider for example
$$F(x, y, z) = (x + y + z)(-x + y + z).$$
The matrix $A$ equals
$$A=\left(\begin{array}{rrr} 1 & 1 & 1 \\ -1 & 1 & 1\end{array}\right),$$
and thus $\gcd(A) = 2$. Although each of the linear forms in the product represents all integers,
it is easy to check that, for instance, $F(x, y, z)$ does not represent $6$.
\medskip

\noindent
{\bf Acknowledgement:} We thank Levent Alpoge and the anonymous referees for some very helpful remarks.

%
%


\begin{thebibliography}{99}



\bibitem{borosh}
I. Borosh, M. Flahive, D. Rubin, and B. Treybig, {\em A sharp bound for solutions of linear
Diophantine equations}. Proc. Amer. Math. Soc. 105(4), 844--846 (1989).

\bibitem{bde}
T.~D. Browning, R.~Dietmann, and P.~D.~T.~A. Elliott,
{\em Least zero of a cubic form}.
 Math. Ann. 352(3), 745--778 (2012).

\bibitem{cassels}
J.~W.~S. Cassels,
 {\em An Introduction to the Geometry of Numbers}.
 Corrected reprint of the 1971 edition, Classics in Mathematics,
  Springer-Verlag, Berlin, 1997.

\bibitem{dietmann}
R.~Dietmann,
 {\em Small solutions of quadratic {D}iophantine equations}.
 Proc. London Math. Soc. 86(3), 545--582 (2003).

\bibitem{cassels_overview}
L.~Fukshansky,
 {\em Heights and quadratic forms: on {C}assels' theorem and its
  generalizations}.
 In W.~K. Chan, L.~Fukshansky, R.~Schulze-Pillot, and J.~D. Vaaler,
  editors, Diophantine methods, lattices, and arithmetic theory of
  quadratic forms, Contemp. Math., Vol. 587, pages 77--94. Amer. Math. Soc.,
  Providence, RI, 2013.


\bibitem{heger}
I. J. Heger, {\em \"Uber die Aufl\"osung eines Systems von mehreren unbestimmten Glei\-chungen des ersten Grades in ganzen Zahlen,
welche eine gr\"ossere Anzahl von Unbekannten in sich schliessen, als sie zu bestimmen verm\"ogen}.
Sitzungsber. Akad. Wiss. Wien (Math.) 21, 550--560 (1856).

\bibitem{jp_jones}
J.~P. Jones,
 {\em Undecidable {D}iophantine equations}.
 Bull. Amer. Math. Soc. (N.S.) 3(2), 859--862 (1980).

\bibitem{masser}
D.W.~Masser,
 {\em Search bounds for {D}iophantine equations}.
 A panorama of number theory or the view from Baker's garden (Zürich, 1999), 247–259, Cambridge Univ. Press, Cambridge, 2002.

\bibitem{mat}
Yu.~V. Matiyasevich,
 {\em The {D}iophantineness of enumerable sets}.
 Soviet Math. Dokl. 11,  354--358 (1970).

\bibitem{zurich}
G.~Maze, J.~Rosenthal, and U.~Wagner,
 {\em Natural density of rectangular unimodular integer matrices}.
 Linear Algebra Appl. 434, 1319--1324 (2011).


\bibitem{oleary}
R. O'Leary and J. D. Vaaler, {\em Small solutions to inhomogeneous linear equations over number
fields}. Trans. Amer. Math. Soc. 336(2), 915--931 (1993).

\bibitem{zhan}
X.~Zhan,
 {\em Completion of a partial integral matrix to a unimodular matrix}.
 Linear Algebra Appl. 414, 373--377 (2006).

\end{thebibliography}


\end{document}